\magnification=\magstep1
\let\iteM=\item
\input amstex.tex

\def\makeatletter{\catcode`\@=11}
\def\makeatother{\catcode`\@=12}
\long\def\comment#1\endcomment{}
\def\lable#1{{\rom [#1]}}
\def\nolables{\def\lable##1{}}
\newcount\sectioncount
\newcount\commoncount\commoncount=1
\newcount\lemmacount\lemmacount=1
\newcount\defcount\defcount=1
\newcount\eqcount\eqcount=1
\newcount\theoremcount\theoremcount=1
\newcount\Ccount\Ccount=1

\def\firstsection#1{\sectioncount=#1 \advance\sectioncount by -1}
\firstsection1
\def\renewcount{\commoncount=1\eqcount=1
}
\renewcount

\def\newsection#1#2\par{\global\advance \sectioncount by 1%
\renewcount%
\specialhead\bigbf{\the\sectioncount}.\lable{#1}\bigbf\ #2
\expandafter\xdef\csname Section #1 \endcsname{%
\the\sectioncount}%
\endspecialhead
\relax}
\def\section#1{{\sl Section \csname Section #1 \endcsname}}

\def\newgrphno#1{
\the\sectioncount.\the\commoncount
\expandafter\xdef\csname Grph #1 \endcsname{%
\the\sectioncount.\the\commoncount}%
\global\advance\commoncount by1\relax}
\def\grph#1{{\csname Grph #1 \endcsname}}

\def\newtheorem#1{
\proclaim{\the\sectioncount.\the\commoncount. Theorem}\lable{#1}
\expandafter\xdef\csname Theorem #1 \endcsname{\the\sectioncount%
\the\commoncount}%
\global\advance\commoncount by1\it\relax}
\def\theorem#1{{\sl Theorem \csname Theorem #1 \endcsname}}

%

%

\def\newlemma#1{
\proclaim{\the\sectioncount.\the\commoncount.\lable{#1} Lemma}
\expandafter\xdef\csname Lemma #1 \endcsname{%
\the\sectioncount.\the\commoncount}%
\global\advance\commoncount by1\it\relax}
\def\lemma#1{{\sl Lemma \csname Lemma #1 \endcsname}}
\def\lemmano#1{\csname Lemma #1 \endcsname}

\def\newproposition#1{
\proclaim{\the\sectioncount.\the\commoncount.\lable{#1} Proposition}
\expandafter\xdef\csname Proposition #1 \endcsname{%
\the\sectioncount.\the\commoncount}%
\global\advance\commoncount by1\it\relax}
\def\proposition#1{{\sl Proposition \csname Proposition #1 \endcsname}}
\def\propositionno#1{\csname Proposition #1 \endcsname}

\def\newcorollary#1{
\proclaim{\the\sectioncount.\the\commoncount.\lable{#1} Corollary}
\expandafter\xdef\csname Corollary #1 \endcsname{%
\the\sectioncount.\the\commoncount}%
\global\advance\commoncount by1\it\relax}
\def\corollary#1{{\sl Corollary \csname Corollary #1 \endcsname}}
\def\corollaryno#1{\csname Corollary #1 \endcsname}

\def\newdefinition#1{
{\smallskip\noindent\bf\the\sectioncount.\the\commoncount.\lable{#1}
Definition.}
\expandafter\xdef\csname Definition #1 \endcsname{%
\the\sectioncount.\the\commoncount}%
\global\advance\commoncount by1\relax}
\def\definition#1{{\sl Definition \csname Definition #1 \endcsname}}
\def\definitionno#1{\csname Definition #1 \endcsname}
\let\dfntn=\definition

\def\makeeq#1{
\expandafter\xdef\csname Equation #1 \endcsname{%
\the\sectioncount.\the\eqcount}
\expandafter\xdef\csname Equationlable #1 \endcsname{{#1}}
\global\advance\eqcount by1\relax}
\def\eqqno#1{\csname Equation #1 \endcsname}

\def\eqlable#1{{\rom{[\csname Equationlable #1 \endcsname]}}}

\documentstyle{amsppt}
\let\item=\iteM 
{\catcode`\@=11\global\def\output@{\shipout\vbox{%
 \iffirstpage@ \global\firstpage@false
  \pagebody 
 \else \ifrunheads@ \makeheadline \pagebody
       \else \pagebody \makefootline \fi
 \fi}%
 \advancepageno \ifnum\outputpenalty>-\@MM\else\dosupereject\fi}
}

\def\example#1{{\sl Example \csname Example #1 \endcsname}}
\def\qed{\ \ \hfill\hbox{$\square$}\smallskip}

\let\definition=\dfntn

\def\rom{\ifmmode \fam\rmfam \else\rm \fi}

\def\longpoints{\leaders\hbox to 0.5em{\hss.\hss}\hfill \hskip0pt}

\input epsf.tex
\newdimen\xsize
\newdimen\oldbaselineskip
\newdimen\oldlineskiplimit
\xsize=.6\hsize
\makeatletter
\def\xbig[#1]#2{{\hbox{$\m@th\left#2\vbox to#1\xsize{}%
\right.\n@space$}}}
\makeatother
\def\xlar[#1]#2{%
\smash{\mathop{ \hbox to #1\xsize{\leftarrowfill}}\limits^{#2}}}
\def\xrar[#1]#2{%
\smash{\mathop{ \hbox to #1\xsize{\rightarrowfill}}\limits^{#2}}}
\def\xline[#1]{\hbox to #1\xsize{\leaders\hrule\hfill}}

\def\restorelineskip{\baselineskip=\oldbaselineskip%
\lineskiplimit=\oldlineskiplimit}
\def\putm[#1][#2]#3{
\hbox{\vbox to 0pt{\parindent=0pt%
\vskip#2\xsize\hbox to0pt{\hskip#1\xsize $#3$\hss}\vss}}}%
%
%
\def\putt[#1][#2]#3{
\vbox to 0pt{\noindent\hskip#1\xsize\lower#2\xsize%
\vtop{\restorelineskip#3}\vss}}

\def\Month{\ifcase\month \or January\or February\or March\or April\or May\or
June\or July\or August\or September\or October\or November\or December\fi}

\font\bigbf=cmbx10 scaled 1200

\thinmuskip = 2mu
\medmuskip = 2.5mu plus 1.5mu minus 2.1mu  
\thickmuskip = 4mu plus 6mu
\font\teneusm=eusm10
\font\seveneusm=eusm7
\font\fiveeusm=eusm5
\newfam\eusmfam
\textfont\eusmfam=\teneusm
\scriptfont\eusmfam=\seveneusm
\scriptscriptfont\eusmfam=\fiveeusm
\def\scr#1{{\fam\eusmfam\relax#1}}
\font\tenmib=cmmib10
\font\sevenmib=cmmib7
\font\fivemib=cmmib5
\newfam\mibfam
\textfont\mibfam=\tenmib
\scriptfont\mibfam=\sevenmib
\scriptscriptfont\mibfam=\fivemib

\font\tensf=cmss10
\font\sevensf=cmss10 scaled 833
\font\fivesf=cmss10 scaled 694
\newfam\sffam
\textfont\sffam=\tensf
\scriptfont\sffam=\sevensf
\scriptscriptfont\sffam=\fivesf
\def\sf{\fam\sffam}
\font\sevensl=cmsl10 scaled 833
\font\fivesl=cmsl10 scaled 694
\scriptfont\slfam=\sevensl
\scriptscriptfont\slfam=\fivesl
\def\sl{\fam\slfam\tensl}
\font\mathnine=cmmi9
\font\rmnine=cmr9
\font\cmsynine=cmsy9
\font\cmexnine=cmex10 scaled 913
\font\ninesl=cmsl9
\font\ninesf=cmss9
\font\ninemib=cmmib9
\def\msmall#1{\hbox{$\displaystyle
\textfont0=\rmnine   \textfont1=\mathnine
\textfont2=\cmsynine \textfont3=\cmexnine
\textfont\slfam=\ninesl \textfont\sffam=\ninesf
\textfont\mibfam=\ninemib{#1}$}}

\hyphenation{Lip-schit-zian Lip-schitz}

\def\cc{{\Bbb C}}

\def\rr{{\Bbb R}}

\def\pp{{\Bbb P}}
\def\ss{{\Bbb S}}

\def\st{_{\rom st}}

\def\exp{{\sf exp}}

\def\id{{\sf Id}}

\def\lim{\mathop{\sf lim}}

\def\sup{{\sf sup}\,}

\def\1{\bold{1}}
\def\eps{\varepsilon}

\def\d{\partial}
\def\barr#1{\mskip1mu\overline{\mskip-1mu{#1}\mskip-1mu}\mskip1mu}

\def\ddef{\mathrel{{=}\raise0.23pt\hbox{\rm:}}}
\def\deff{\mathrel{\raise0.23pt\hbox{\rm:}{=}}}
\def\ge{\geqslant}

\def\<{\langle}
\def\>{\rangle}
\def\fraction#1/#2{\mathchoice{{\msmall{ #1\over#2}}}%
{{ #1\over #2 }}{{#1/#2}}{{#1/#2}}}

\def\le{\leqslant}


\def\mapdown#1|#2{\llap{$\vcenter{\hbox{$\scriptstyle #1$}}$}
 {{ \big\downarrow}}
  \rlap{$\vcenter{\hbox{$\scriptstyle #2$}}$}}

\def\emptyset{\varnothing}

\def\state#1. {\smallskip\noindent{\bf#1. }}
\def\qed{\ \hbox{ }\ \hbox{ }\ {\hfill$\square$}}
\def\Chi{\raise 2pt\hbox{$\chi$}}
\let\phI=\phi\let\phi=\varphi\let\varphi=\phI

%
%
%
%
%
%
%
%
%
%
%
%
%
 %
%
%
%
%
%
%
%
%
%
%
%
%
%
%
%
%
%

\def\ie{{\sl i.e.\/}\ \hskip1pt plus1pt}

\def\.{\thinspace}
\def\3{\ss}

\def\sli{{\sl i)} }             
\def\slii{{\sl i$\!$i)} }       
\def\sliii{{\sl i$\!$i$\!$i)} }

\def\cal#1{{\scr{#1}}}
\def\cala{{\cal A}}
\def\alg(#1){\cala(\barr{#1})}

\def\calm{{\cal M}}


\voffset=-.27truein 

\voffset=-.2in
\vsize=580pt
\baselineskip=12.33pt plus .2pt

\nolables

\long\def\comment#1\endcomment{}

\medskip
\centerline{\bigbf COMPLETE HYPERBOLIC NEIGHBORHOODS}
\medskip
\centerline{\bigbf IN ALMOST-COMPLEX SURFACES}

\medskip
\centerline{\rm R. Debalme$^*$, S. Ivashkovich$^{*,\dag}$}

\thanks
$^*$ Address:
 UFR de Math\'ematiques, Universit\'e de Lille 1, 59655 Villeneuve d'Ascq Cedex, 
 France. e-mail: debalme\@agat.univ-lille1.fr; ivachkov\@agat.univ-lille1.fr 
\newline\indent
$^\dag$ This work was partially done when the author visited the 
university of Wisconsin-Madison. He would like to thank this Institution for the 
hospitality.
\endthanks

\endtopmatter

\def\version{\hbox{\font\fiverm=cmr5 \fiverm
Version of \the\day.\the\month.}}
\leftheadtext{\hss\vtop{%
\line{\hfil R.\.Debalme\ \ S.\.Ivashkovich \hfil
\llap{\version}}%
\vskip 4pt \hrule }\hss}

\rightheadtext{\hss\vtop{%
\line{\rlap{\version}%
\hfil Complete hyperbolic\ neighborhoods\hfil}%
\vskip 4pt \hrule }\hss}

\subhead 0. Introduction \endsubhead

\smallskip

Our goal here is to prove that each point in an almost-complex surface 
has a basis of complete  hyperbolic neighborhoods. The problem is local, and 
therefore we can consider the case when our surface is $\Bbb R^4$ with an 
arbitrary almost-complex structure $J$. Throughout this paper, almost-complex 
structures are assumed to be of class $C^{1,\alpha }$ for some 
$0<\alpha <1$. Let $C$ be some non-singular $J$-complex curve passing through 
the origin.

We shall prove the following

\state Theorem. {\it There exists a basis $\{ U_j\} $ of neighborhoods of 
zero in $\rr^4$, such that:

\smallskip
1) $(U_j,J)$ are complete hyperbolic in the sence of Kobayashi;

\smallskip
2) $(U_j\setminus C,J)$ are complete hyperbolic as well.
}

\smallskip In the case when $J$ is integrable, this result follows easily from 
the Pick-Schwarz lemma. Indeed, one can choose a local holomorphic  coordinate 
system so that $C$ is one of the axes, and then set $U_j=\Delta^2_{{1\over j}}$,
the bidisk of radius ${1 \over j}$. Now the Kobayashi distance  on $U_j$ and 
$U_j\setminus C$ 
can be 
calculated explicitly and one can see from this explicit form that it is 
complete, see e.g. [Ko-2], Ch.2. The utility of complete hyperbolic neighborhoods 
in the theory of complex hyperbolic manifolds, first observed by Kiernan and 
Kobayashi, has become standard since. The fact that this result 
remains true for {\it any} almost-complex structure is somewhat suprising. 
Really, given {\it any} germ of a non-singular real surface $C\ni 0$ in 
$\rr^4$, one can easily construct an almost-complex structure $J$ in a 
neighborhood of zero such that $C$ becomes a $J$-complex curve. Note also 
that because of the non-existence of $J$-holomorphic functions for general
$J$, the hyperbolicity of $U\setminus C$ for a $J$-complex curve 
$C$ is not evident either.   

\smallskip An open subset $Y$ in an almost-comlex manifold $(X,J)$ is called 
{\sl locally complete hyperbilic} (l.c.h.) if for every $y\in \bar Y$ there 
exists a neighborhood $V_y\ni y$ such that $V_y\cap Y$ is complete hyperbolic.
An example of this situation we shall have in mind is when $(X,J)$ is an 
almost-complex surface and $Y=X\setminus D$, where $D=\cup_kD_k$ is a 
(reducible) 
curve with irreducible $J$-complex components $D_k$, which do not have 
cusps. Our theorem insures that such a $Y$ is l.c.h. 

The following proposition is due to Zaidenberg  [Za-1] in the case of complex 
manifolds . Its proof goes through  for  almost-complex
manifolds as well.

\break

\state Proposition. {\it  Let $Y$ be a relatively compact l.c.h. open 
subset of an almost-complex manifold $(X,J)$. $Y$ is hyperbolically 
imbedded into $X$  if and only if $Y$ does not contain  $J$-complex 
lines and admits no limiting $J$-complex lines. 
}

Recall that an open subset $Y$ of an almost-comlex manifold $X$ is called 
{\sl hyperbolically imbedded} into $X$ if for any two sequences $\{ x_n\} , 
\{ y_n\} $ in $Y$ converging to $x\in \bar Y$ and   $y\in \bar Y$, 
respectively, one has that ${\overline \lim_{n\to \infty }}k_{Y,J}(x_n,y_n)>0$. Here 
$k_{Y,J}$ denotes the Kobayashi pseudodistance on the manifold $(Y,J)$.
It it worth observing that if $Y$ is hyperbolically imbedded 
into $(X,J)$ and is l.c.h. then $(Y,J)$ is complete hyperbolic, see [Ki]. As 
in the integrable case, by a $J$-complex line we understand 
a non-constant $J$-holomorphic map $f:\cc \to X$ such that $\Vert d_zf({\d 
\over \d x})\Vert_h\le 1$ for all $z\in \cc $. Here $h$ stands for some 
$J$-Hermitian metric on $X$. A $J$-complex line $f$ is said to be a {\sl 
limiting line} for $Y$ if $f(\cc )\subset \d Y$ and for every radius $R$ 
there exists a sequence $f_n^R:\Delta (R)\to Y$ of $J$-holomorphic maps 
of the disk of radius $R$ into $Y$ converging uniformly to $f\mid_{\Delta (R)
}$. 

Using the Brody reparametrization lemma we can (and we 
shall always) assume that $\sup \{ \Vert d_zf({\d\over \d x})\Vert_h: 
z\in \cc \} = $ $\Vert d_0f({\d \over \d x})\Vert_h =1$.

\smallskip
An immediate consequence of our theorem and this proposition is the following

\state Corollary 1. {\it Let $(X,J)$ be a compact almost-complex manifold of 
dimension $4$. Let $D= \displaystyle \cup_{j=1}^n D_j$ be a reducible $J$-
complex curve, each 
irreducible component $D_j$ of which is an immersed $J$-complex curve. Suppose 
that for every $j=1, \ldots ,n$ the  curve $(D_j \setminus Sing(D),J\mid_{D_j})$ is 
hyperbolic.
Then  $(X \backslash D,J$) is hyperbolically imbedded into $(X,J)$ if and only 
if there is no 
$J$-complex line in $ X \backslash D$. 
}

\smallskip By an immersed $J$-complex curve we mean here the image of 
a $J$-holo\-morphic map $u:S\to X$ of a compact Riemann surface $S$, 
such that $du$ never vanishes (\ie, $u(S)$ is a curve without cusps).
Let us state one more corollary, which in fact was our initial motivation for 
considering the existence of complete hyperbolic neighborhoods.

\state Corollary 2. {\it Let $\calm_{\omega ,5l}$ be the Banach manifold 
consisting of pairs $(J,\{ D_j\}_{j=1}^5)$, where $J$ is any almost-complex 
structure on $\cc\pp^2$ tamed by the Fubini-Studi form $\omega $ and 
$\{ D_j\}_{j=1}^5 $ the union of five $J$-complex 
lines in $\cc\pp^2$ in general position. The set ${\cal H}_{\omega ,5l}$ 
  consisting of $(J, \{ D_j\}_{j=1}^5)$ with $Y=(\cc\pp^2\setminus 
 \bigcup_{j=1}^5 D_j,J)$ hyperbolically imbedded into $(\cc\pp^2, J)$
 is an open nonempty subset of $\calm_{\omega ,5l}$. 
}

Nonemptyness of ${\cal H}_{\omega ,5l}$ is the theorem of Bloch. Really, 
 for the standard  complex structure on $\cc\pp^2$ and standard complex 
lines $D_1,...,D_5$ in general position the complex manifold $\cc\pp^2\setminus 
D$ is Kobayashi hyperbolic. Therefore we obtain plenty of 
examples of almost-complex hyperbolic manifolds of this type (\ie compact 
manifold minus a reducible $J$-complex curve).

\bigskip\bigskip\bigskip\break

\smallskip We are interested here in the case of almost-complex surfaces, 
because a generic almost-complex structure $J$ on a manifold of higher dimension 
does not possess any $J$-complex hypersurface even locally.

\smallskip We would like also to point out that the following convention is 
applied throughout this paper (and was already used in this introduction): 
if we use some statement or definition and refer to the original paper 
(where this statement/definition is proved/formulated for the integrable 
case), this means that the same proof/defini\-tion goes through also in 
the non-integrable case and no additional discussion is needed. The 
reader is supposed to look at the original paper or, if 
she/he prefers to have everything in one place the book [Ko-2] 

\smallskip The composition of the paper is the following. In \S 1 we recall 
very briefly a few facts from (almost)-complex hyperbolic geometry. In 
\S 2 we state the principal result of this paper, a sort of "gauge invariant 
Schwarz lemma", and deduce  theorem from this lemma. In 
\S 3 we give the proof of this Schwarz lemma. In the last paragraph we give 
the proofs of  Corollaries 1 and 2 and some others results.

\bigskip
\centerline{\bigbf Table of Contents}

\bigskip\bf
0. Introduction.

\smallskip
1. Kobayashi hyperbolicity of almost-comlex manifolds.

\smallskip
2. Local complete hyperbolicity in $(\rr^4 , J)$.

\smallskip
3. Proof of the Schwarz-type lemma.

\smallskip
4. Corollaries.

\medskip
References

\bigskip\noindent\bf 
1. Kobayashi hyperbolicity of almost-complex manifolds.

\smallskip\rm
Let $(X,J)$ be an almost-complex manifold. For every two sufficiently 
close points $p$ and $q$ in $X$, there exsists a $J$-holomorphic mapping 
$u:\Delta \to X$ such that $u(0)=p$ and $u(a)=q$, where $a\in ]0,1[$ (see 
[De], Lemma 1). This allows one to define the Kobayashi pseudodistance $k_{X,J}$ 
on $(X,J)$. 

For every point $p\in X$ and every tangent vector $v\in T_pX$, there exists 
for some $R>0$ a 
(non-unique) $J$-holomorphic map $u:\Delta (R)\to X$ with $u(0)=p$ and 
$d_0u({\d \over \d x})=v$. It is therefore possible to introduce the Royden 
infinitesimal 
pseudonorm and prove that the integral pseudometric generated 
by this pseudonorm coincides with the Kobayashi pseudometric. In the case when 
this pseudometric is actually a metric (\ie $(X,J)$ is Kobayashi hyperbolic), 
this metric induces the given topology on $X$ (by a theorem of Barth [Bt], 
see also [Kr]).

If $(X,J)$ is not hyperbolic then the result of Brody [Br] states  that
  
\bigskip\bigskip\bigskip\break

$$
\sup \{ \Vert d_0f(\d / \d x)\Vert_h  :  f \in {\Cal O}_J (\Delta ; X)\} 
= \infty ,\eqno(1.1)
$$ 
where ${\Cal O}_J(\Delta ; X)$ denotes the space of $J$-holomorphic maps 
from the unit disk $\Delta $ to $X$, and $h$ is a $J$-Hermitian metric 
on $X$.

The proof of the following fundamental result, known as Brody's 
reparametrization lemma, does not 
uses the integrability of $J$. Let $(X,J)$ be an 
almost-complex manifold. 
Let $f : \Delta (R) \to X$ be a $J$-holomorphic map  with 
$\Vert d_0f({\d \over \d x})\Vert  \geq c > 0$.
Then there exists a $J$-holomorphic map $\tilde f :\Delta (R) \to X$
such that

$$
\displaystyle\sup_{z \in \Delta (R)} \Vert d_z\tilde f(\d / \d x)\Vert_h 
\bigl[ {R^2-\vert z\vert ^2 \over R^2}\bigr] =
\bigl|| d_0\tilde f(\d / \d x)\Vert  = c. \eqno(1.2)
$$

\smallskip A combination of (1.1) and (1.2) implies that a compact 
almost-complex manifold $(X,J)$ is not Kobayashi 
hyperbolic if and only if it contains a $J$-complex line, see [De] or [KrOv] 
for more details.

\bigskip

\subhead 2. Local complete hyperbolicity in $(\Bbb R^4,J)$ \endsubhead

\smallskip

\smallskip We are going to construct a basis of complete hyperbolic 
neighborhoods of the origin in $(\rr^4,J)$, where $J$ is an arbitrary almost-complex 
structure of class $C^{1,\alpha }$, $0<\alpha <1$.

There exists a local coordinate system 
$(z_1,z_2)$ in a neighbourhood of $0$ such that our almost-complex 
structure in these coordinates has the form 

$$
J(z_1,z_2) = \left( \matrix A(z_1,z_2)&0 \\ 0&B(z_1,z_2) \\ \endmatrix \right),
$$

\noindent
and $J(0)=J_{st }$, the standard structure of $\Bbb C^2$. If moreover, a 
non-singular $J$-complex curve $C\ni 0$ is given, we can arrange that in this 
coordinate system $C$ is defined by $z_1=0$, see [Sk].

 Consider the structure $J_{\varepsilon }(z_1,z_2):=J(\varepsilon z_1,\varepsilon
z_2)$. All we need to prove is that, for $\varepsilon $ small enough, the unit 
bidisk $\Delta^2:=\{ (z_1,z_2): \vert z_1\vert <1, \vert z_2\vert <1\} $ and 
$\Delta^2\setminus \{ z_1=0\} $ are complete 
hyperbolic with respect to $J_{\varepsilon }$.
Denote by $z=x+iy$ the coordinate on $\Delta$.  
Let $f=(u,v) \in {\Cal O}_J(\Delta,\Delta^2)$. The Cauchy-Riemann equation for 
the first component is 

$$
{\partial u \over \partial x} + A(u(z),v(z)){\partial u \over \partial y}=0,
$$
where $A(u(z),v(z))$ is a 2x2 real matrix. An obvious transformation gives


$${\partial u \over \partial \bar{z}}+(1-Ai)^{-1}(1+ A i){\partial u \over 
\partial z}=0.
$$

Denote $q_A(z_1,z_2):=-[1-A(z_1,z_2)\circ J_{st}]^{-1}[1 + A(z_1,z_2)\circ 
J_{st} ] $.
For $z,u,v$ fixed, $q_A(z,v,u)$ is an $\rr $-linear operator. Let us decompose 
it into   
$\Bbb C$-linear and $\Bbb C$-antilinear part $q_A(z,v,u)=\mu^1(z,v,u) + 
\mu^2(z,v,u) \circ \sigma$, where $\sigma$ is the conjugation operator.
Then the Cauchy-Riemann equation for the first component becomes

\bigskip\bigskip\bigskip\break

$$
{\partial u \over \partial \bar{z}}-\mu^1(z,v,u){\partial u \over \partial z} -
\mu^2(z,v,u) \overline{ {\partial u \over \partial z}}=0. \eqno(2.1)
$$

The second function $v$ will be considered as a parameter. For functions, tensors
etc. defined on the open subset $W\subset \rr^N$, by the $C^k(W)$-norm of an 
object $T$ 
 we shall mean 
$\Vert T\Vert_{C^k(W)}:=\sum_{r=0}^k\sup \{ \Vert D^r(T(x)\Vert :x\in W\} $. 

The following Schwarz-type 
lemma  due to Gromov gives us a $C^1$-bound on the parameter $v$.

\state Lemma 2.1. {\it If $J$ is sufficiently close to $J\st $ in $C^1$-sense 
on $\Delta^2$, then there exists a constant $C=C(\Vert J-J\st 
\Vert_{C^1(\Delta^2)})$ such that every $J$-holomorphic map 
$f:\Delta \to \Delta^2$ satisfies 

$$
\Vert df(0)\Vert \le C. \eqno(2.2)
$$
}
 
Therefore there exists (another) $C$ such that for every  $J$-holomorphic
$f=(u,v):\Delta \to \Delta^2$ we have that 
$\Vert u\Vert_{C^1(\Delta ({1\over 2}))}$ and 
$\Vert v\Vert_{C^1(\Delta ({1\over 2}))}$ are bounded  by $C$. Rescaling 
by $z\to {1\over 2}z$ we can suppose that this bound holds everywhere 
on $\Delta $. Because of the local character of  our considerations 
we can use the standard Euclidean norm $\Vert \cdot \Vert $ in all 
estimates like (2.2).

\smallskip Suppose a relatively compact domain $D\subset \cc $ is given. 
Fix some universal covering map $\pi :\Delta \to D$. Suppose $\mu^1, \mu^2\in 
C^1(\rr^6)$ are given and have finite (in the sequel sufficiently small) 
$C^1(\rr^6)$-norm. Let us consider the generalized Beltrami equation (2.1)
for $C^1$-mappings $u:\Delta \to D$. We  suppose that the $C^1(\rr^2)$-
norm of the parameter $v$ is bounded by some $C$. Further, for a $C^1$-map 
$u:\Delta \to D$ denote by $u_{\pi }$ any single-valued branch of $\pi^{-1}
\circ u$. Denote by $\lambda $ the coordinate on the universal cover $\Delta $.
Let also $\phi_a(\lambda )={a-\lambda \over 1-\bar a\lambda }$ be the 
authomorphism of the unit disk interchanging $a\in \Delta $ with the origin.  
Equation (2.1) rewrites for $u_{a,\pi }:=\phi_a\circ u_{\pi }$ as 

$$
{\partial u_{a,\pi } \over \partial \bar{z}}-\mu^1_{\pi }(z,v,a,u_{a,\pi })
{\partial u_{a,\pi } \over \partial z} -
\mu^2_{\pi }(z,v,a,u_{a,\pi }) \overline{ ({\partial u_{a,\pi } \over 
\partial z})}=0, \eqno(2.3)
$$
where $\mu^1_{\pi }(z,v,a,u_{a,\pi })=\mu^1(z,v,(\pi \circ \phi_a )(u_{a,\pi }))
$ and $\mu^2_{\pi }(z,v,a,u_{a,\pi })=\mu^2(z,v,(\pi \circ \phi_a )(u_{a,\pi })
\cdot \bigl(\overline{{\d (\pi \circ \phi_a )\over \d \lambda }}\cdot 
\bigl[{\d (\pi \circ \phi_a )\over \d \lambda }\bigr]^{-1}\bigr)$.

In the next paragraph we shall prove the following Schwarz-type lemma:

\bigskip\bigskip\bigskip\break

\state Lemma 2.2. {\it There exists an $\eps = \eps (\pi , C)>0$ and 
a constant $K=K(\pi, C,\eps )<\infty $, such that for every $C^1$-solution $w:\Delta \to \Delta $ 
of (2.3) with coefficients and parameters satisfying 
$\Vert \mu^{1,2}\Vert_{C^1(\rr^6)})<\eps $, $\Vert v\Vert_{C^1(\Delta )}<C$ 
and such that $w(0)=0$, we have the following estimate  

$$
\Vert dw(0)\Vert \le K.\eqno(2.4)
$$.
}

\state Remark. {\bf 1.} {\rm When one proves the statement, which is now called 
the Schwarz-Pick lemma for holomorphic mappings $u:\Delta \to \Delta $, 
one proves in fact two statemets. Firstly, the Schwarz lemma, which gives 
the universal bound for the absolute value of the derivative of $u$ 
provided $u(0)=0$, \ie, one obtains $\vert du(0)\vert \le 1$. Second, using 
"the gauge invariance" of the homogeneous 
 Cauchy-Riemann equation with respect to the "gauge group" 
 $G:=Aut_{Hol}(\Delta )$, 
 one obtains the "gauge invariant" form of this statement: 
 
 $$
 {\vert du(0)\vert \over 1-\vert u(0)\vert^2 }\le 1, \eqno(2.5)
 $$
 which means exactly the  completeness of Kobayashi metric on the unit disk.
 
 The generalized Cauchy-Riemann  (or Beltrami) equation (2.1) has no symmetries
 except for very special $\mu $'s. Therefore, instead of considering only 
 one equation we "symmetrize" the problem and consider the whole family of 
 equations (2.3) with a parameter $\phi_a\in G$, our gage group. Uniform 
 estimate (2.4) of the derivative at zero allows us to make the 
 same conclusions as for the standard Cauchy-Riemann equation, because the 
 problem is now "enforced" with nessessary symmetries.
 
 \smallskip\noindent\bf 2. \rm 
 However the price for this is that the coefficients of (2.4) are no longer 
 bounded in any Sobolev or H\"older space as the parameter $a\to \d \Delta $.
 They are  bounded in the $\sup $-norm only, and this is clearly not sufficient 
 to bound the derivative at zero. The trick is (see Step 1 of the proof in the 
 next 
 paragraph) that there is an {\it a priori} bound for the solutions $w$ of (2.4) 
 satisfying $w(0)=0$ in the H\"older norm on a neighborhood of zero. This  
 allows us to control the H\"older norm of the coefficients and apply  
 elliptic bootstrap one more time (see Step 2).
 }     

\smallskip
\state Corollary 2.3. {\it For $J$ sufficiently $C^1$-close to $J\st $, the 
open set $(\Delta^2,J)$ is complete hyperbolic.
}
\state Proof. {\rm As we have already explained, by taking $J$ close to $J\st $ 
we can bound the components $(u,v)$ of a $J$-holomorphic map 
$f:\Delta \to \Delta^2$ using Gromov's lemma. The explicit form of $q_A$ and 
consequently of $\mu^{1,2}$ shows that, again by taking $J$ close to $J\st $, 
we can ensure that $\Vert \mu^{1,2}\Vert_{C^1}<\varepsilon $ for the $\eps $ 
from Lemma 
2.2 applied to the identity covering $\pi =\id :\Delta \to \Delta $. Let 
$u(0)=a$. Applying Lemma 2.2 to $w=\phi_a\circ u$, we obtain that 
$\Vert dw(0)\Vert 
\leq K$, where $K$ does not depend on $w$.
 
\bigskip\bigskip\bigskip\break

But $dw(0)={1\over \vert a\vert^2 - 1
}\cdot du(0)$ and therefore, for every solution $u : \Delta \mapsto \Delta$  
of $(2.1)$ such that $u(0)=a$, one has  ${1 \over ||du(0)||} \geq {1 \over K 
(1-|a|^2)}$. This means that, for every 
vector $(\xi, 0)$ in $T_{(a,b)}\Delta^2$, we have $k_{\Delta^2,J}((a,b),(\xi,
0)) \geq {C_1 |\xi|\over 1-|a|^2}$, by the definition of the Kobayashi-Royden 
pseudometric.

The same considerations apply, of course, to the second component of $(u,v)$. 
This implies that   
$ k_{\Delta^2,J}((a,b),(0 ,\eta)) \geq {C_2 |\eta|\over 1-|b|^2}$. Therefore,

$$ 
k_{\Delta^2,J}((a,b),(\xi,\eta)) \geq \max({C_1 |\xi|\over 1-|a|^2},{C_2 |\eta|
\over 1-|b|^2}), \eqno(2.6)
$$
which gives us the desired conclusion.
}
\smallskip
\hfill 
\qed

\state Corollary 2.4. {\it For $J$ sufficiently $C^1$-close to $J\st $ on the 
whole bidisk 
$\Delta^2$, the open set $(\Delta^2\setminus (\{ 0\} \times \Delta ), J)$ is 
complete 
hyperbolic.
}
\state Proof. {\rm 
Let us apply Lemma 2.2 to the universal covering $\pi :\Delta \to \Delta^*=D$ 
given by $\pi (\lambda ) =\exp{{\lambda -1\over \lambda +1}}$.  
For a solution $u:\Delta \to \Delta^*$ of (2.1), set $a=u(0)$ and $b=ln(a)$ with 
imaginary part $Im(b)\in [-\pi ,\pi [$. Let also $c={b+1\over 1-b}\in \Delta $.
Therefore $c$ is one of the  preimages of $a$ under the covering map $\pi $. 
Take the single-valued branch $u_{\pi }$ of $\pi^{-1}\circ u$ such that $u_{\pi }(0)=
c$, and consider $w=\phi_c\circ u_{\pi }$, a solution of (2.3). Lemma 2.2 
gives us that $\Vert dw(0)\Vert \le K$, where $K$ does not depend on anything 
involved. But direct computation shows that 

$$
dw(0)= {|c|^2-1 \over (1-|c|^2)^2}{2\over (1+b)^2}{du(0) \over a}
={1 \over |c|^2-1}\cdot {2\over (1+b)^2}{du(0) \over a}.
$$

\noindent
Therefore

$$
dw(0)={1 \over \vert {1+\ln a \over 1-\ln a}\vert^2-1}{2\over (1+\ln a)^2}{du(0) 
\over a},
$$
and hence 

$$
{1 \over |du(0)|} \geq K_1
{ 1\over |a|\ln{1\over |a|}}
$$
for $a$ close to zero.
By the definition of the Kobayashi-Royden pseudometric, for every 
vector $(\xi, 0)$ in $T_{(a,b)}(\Delta^2\setminus \Delta^*)$, we have

$$
k_{\Delta^2\setminus \Delta^*,J}((a,b),(\xi,0)) \geq K_1
{ 1\over |a|\ln{1\over |a|}}\eqno(2.7)
$$
for $a$ close to zero, and therefore, 
$k_{\Delta^2\setminus \Delta^*,J}((a,b),(\xi,\eta)) = 
\bigcirc ({ 1\over |a|\ln{1\over |a|}})\sqrt{\xi^2+\eta^2}$. The consequence is that every path leading to $0$ has 
infinite length, which proves the corollary.

\hfill $\square$

\bigskip\bigskip\bigskip\break

\medskip\noindent\bf
\S 3. Proof of the Schwarz-type lemma.

\smallskip\rm

\state Proof. \rm Suppose that the conclusion of the lemma is not true. 
Let $(h_n)$ be a sequence of mappings $h_n : \Delta \mapsto \Delta,
h_n(0)=0$, solutions of $(2.3)$, such that $\vert dh_n(0)\vert \to \infty $. 
For each $n$, denote by $v_n$ the corresponding functional parameter and by 
$a_n $ the corresponding parameter in $\Delta $. 

\smallskip\noindent\sl
Step 1. \it  There exists $\delta >0$ and $n_0$ such that, for $n\ge n_0$ 
$h_n\in C^{0,{1\over 2}}(\Delta_{\delta },\Delta_{{1\over 2}})$ and 
$\{ h_n\}_{n\ge n_0} $ is a  bounded set in  this space.

\smallskip\rm 
Consider  a test function $\rho$ which is $C^{\infty}$ and such that
$\rho \equiv 1$ on $\Delta_{3\over 4}$ and $\rho \equiv 0$ outside $\Delta$, 
and denote
$h_n^{\rho}=\rho h_n$. The equation $(2.3)$ for $h_n$, written in terms of 
$h_n^{\rho}$, gives

$$
{\partial  h_n^{\rho} \over \partial \bar{z}}
-\mu^1_{v,a_n}(z, h_n){\partial  h_n^{\rho} \over \partial z}-
\mu^2_{v,a_n}(z, h_n)
\overline{({\partial  h_n^{\rho} \over \partial z})}=
h_n{\partial  \rho \over \partial \bar{z}}
-h_n\mu^1_{v,a}(z,h_n){\partial \rho \over \partial z}-
$$

$$
-\bar h_n\mu^2_{v,a}(z, h_n)
\overline{({\partial  \rho \over \partial z})}. \eqno(3.1) 
$$
Consider the Cauchy-Green operator $T_{CG} : L^4(\Delta) \to L^{1,4}(\Delta)$ 
and the Calderon-Zygmund operator 
 $T_{CZ} : L^4(\Delta) \to L^4(\Delta)$, defined respectively as
   
$$
T_{CG}g(z) := {1 \over 2i\pi}\int\int_{\Delta}{g(\zeta) \over \zeta - z}d\zeta\wedge 
d\bar{\zeta},
$$
and

$$
T_{CZ} g(z):=
 p.v. {1\over 2i\pi}\int\int_{\Delta}{g(\zeta) \over (\zeta - z)^2}d\zeta\wedge 
d\bar{\zeta}.
$$ 

$T_{CG}$ and $T_{CZ} $ are bounded linear operators in Sobolev and 
H\"older spaces. Let us denote by $C_{k,p}$ the norm of 
$T_{CG}:L^{k,p}(\Delta) \to L^{k+1,p}(\Delta)$ and by $C_{k,\alpha }$ the norm of 
$T_{CG}:C^{k,\alpha }\to C^{k+1,\alpha }$. By $A_{k,p}$ and $A_{k,\alpha }$ we denote 
the similar constants for $T_{CZ} : L^{k,p}\to L^{k,p}$ and $T_{CZ} : C^{k,\alpha }\to 
C^{k,\alpha }$. Here always $k\ge 0, 1<p<\infty $ and $ 0<\alpha <1$.

We know that $\bar{\partial} \circ T_{CG} = Id$ and that $T_{CZ} = \partial \circ T_{CG}$.
We also know that  on the closure in $L^{1,2}$ of the space of compactly supported 
$C^{\infty}$ maps, $T_{CZ} \circ {\bar \partial}= \partial$ and $T_{CG} \circ \bar{\partial}= Id$.
We remark that the $h_n^{\rho}$ are in this space.
Set $g_1={\partial  \rho \over \partial \bar{z}}
-\mu^1_{v,a}(z,h_n){\partial \rho \over \partial z}$ and $g_2=-\mu^2_{v,a}(z, h_n)
\overline{({\partial  \rho \over \partial z})}$.
The equation $(3.1)$ can be written as  
 
$$
[Id -\mu^1_{v,a}(., h_n)\circ T_{CZ} - \mu^2_{v,a}(., h_n)\circ \sigma \circ T_{CZ}]{\partial \over 
\partial {\bar z}} h_n^{\rho}= h_n g_1 + {\overline h_n}g_2.\eqno(3.2) 
$$

We remark that $h_n g_1 + {\overline h_n}g_2$ is a bounded sequence in $L^4(\Delta, \Bbb R^2)$. 
Indeed, $h_n$ takes its values in $\Delta$ and $g_1$ and $g_2$ are bounded.
Since $||T_{CZ} ||_{L^4}$ is bounded, and as $\mu^1$ and $\mu^2$ are as small 
as we want, this operator
is invertible. Therefore,

\bigskip\bigskip\bigskip\break

$$
{\partial \over \partial {\bar z}} h_n^{\rho} =
[Id -\mu^1_{v,a}(., h_n)\circ T_{CZ} - \mu^2_{v,a}(., h_n)\circ \sigma \circ T_{CZ}]^{-1}
( h_n g_1 + {\overline h_n}g_2),
$$
and hence

$$ 
h_n^{\rho} =
T_{CG} \bigl([Id -\mu^1_{v,a}(., h_n)\circ T_{CZ} - \mu^2_{v,a}(., h_n)\circ \sigma \circ T_{CZ}]^{-1}
( h_n g_1 + {\overline h_n}g_2)\bigr). \eqno(3.3)
$$

Thus $(h_n^{\rho})$ is a bounded sequence in $L^{1,4}(\Delta ,\Delta )$.
By Sobolev imbedding  $L^{1,p} \subset C^{0,\alpha ={2\over p}}$, our 
sequence
$(h_n^{\rho})$ is  bounded  in $C^{0,{1 \over 2}}$.
Therefore, $(h_n^{\rho})$ is a bounded and equicontinuous family of mappings on
compacts in $\Delta$. From the theorem of Ascoli we deduce that it is compact, and extract
from it a converging subsequence. Therefore, $(h_n)$ converges uniformly on 
$\Delta_{3\over4}$ to a mapping $h:\Delta_{{3\over 4}} \to \Delta$ such that $h(0)=0$.
There exists $\delta$ such that $h(\Delta_{\delta}) \subset \Delta_{1\over2}$ and 
consequently $h_n(\Delta_{\delta})\subset \Delta_{{1\over 2}}$ for $n$ big enough. 
Therefore  $h_n$ is bounded in  
$C^{0,{1 \over 2}}(\Delta_{\delta},\Delta_{{1\over 2}})$ starting from 
some $n>>1$.

\medskip\noindent\sl
Step 2. \it There exists $n_1$ such that $h_n\in C^{1,{1\over 2}}(\Delta_{\delta 
\over 2})$ for $n\ge n_1$ and $\{ h_n\}_{n\ge n_0}$ is a bounded set in  this 
space.

\smallskip\rm
Parameters $v_n$ are bounded in  $C^1(\Delta)$, and $h_n\in  
C^{0,{1 \over 2}}(\Delta_{\delta},\Delta_{{1\over 2}})$. 
Therefore $\mu^1_{v,a_n}(z,h_n)$ and $\mu^2_{v,a_n}(z, h_n)$ are bounded  in 
$C^{0,{1 \over 2}}(\Delta_{\delta })$. 
Put $f_n(z):=h_n(\delta z)$. This is a bounded sequence in 
$C^{0,{1 \over 2}}(\Delta , \Delta )$. All we need  to prove is that $\{ f_n\} $ 
is bounded in  $C^{1,{1\over 2}}(\Delta_{1 \over 2})$. We see that 
$f_n$ satisfy

$$
{\partial f_n \over \partial \bar{z}}-\mu^1_{v_n,a_n}(z,f_n(z))
{\partial f_n \over \partial z} -
\mu^2_{v_n,a_n}(z,f_n) \overline{ ({\partial f_n \over \partial z})}=0. 
\eqno(3.4)
$$
with $\mu^1_{v_n,a_n}(z,f_n)$ and $\mu^2_{v_n,a_n}(z, f_n)$ uniformly 
bounded in $C^{0,{1\over 2}}(\Delta ,\Delta )$. We can repeat now the reasonings
from Step 1 applied to $f_n$ instead of $h_n$ and arrive at equation 
(3.2) but considered in H\"older spaces. In this way we get the 
boundedness of $f_n$ in $C^{1,{1\over 2}}(\Delta_{{1\over 2}})$. Contradiction.

\hfill $\square$

\medskip\noindent\bf
 4. Corollaries.

\smallskip\rm We shall prove Corollaries 1 and 2 from the Introduction 
first .
The crucial point in the proof is the following statement about the positivity
of intersections of $J$-complex curves due to [M-W]. 

\state Theorem. {\it Let $u_i:\Delta\to(\rr^4, J)$, $i=1, 2$ be two primitive
distinct $J$-complex  disks  such that $u_1(0)=u_2(0)=0$. Set $M_i \deff u_i
(\Delta)$. Let $Q = M_1\cap  M_2$ be the intersection set of the disks.
If $J$ is $C^1$-smooth, then the following is true:

\bigskip\bigskip\bigskip\break

\smallskip
\sli The set $\{\, (z_1,z_2) \in \Delta \times \Delta: 
u_1(z_1)= u_2(z_2)\,\}$ is
a discrete subset of $\Delta\times \Delta $. In particular, $u_1(\Delta) \cap
u_2(\Delta)$ is a countable set. If moreover $u_1(\d \Delta) \cap 
u_2(\d \Delta) = \emptyset $ then  $u_1(\Delta) \cap u_2(\Delta)$ is finite.

\smallskip
\slii The~intersection index $\delta_p$ of $M_1$ and $M_2$ at any such point 
$p\in Q$ is strictly positive.
Moreover, if $\mu_1$ and $\mu_2$ are the~multiplicities of $u_1$ and $u_2$ in
$z_1$ and $z_2$, respectively, with $u_1(z_1)=u_2(z_2)=p$, then the
intersection number $\delta_p$ at $p=u_j(s_j)$ is at least $\mu_1\cdot \mu_2$;

\smallskip
\sliii $\delta_p=1$ iff $M_1$ and $M_2$ intersect at $p$ transversally.
}

For a proof see Theorems 7.1 and 7.3 in [M-W], or Theorem 3.5.1 in [I-S].

We shall need the following topological description of the total intersection number 
$\delta =\sum_{p\in Q}\delta_p$ of $M_1$ and $M_2$.

Take a sphere $S_r$ of radius $r$ around  the origin, which bounds the ball 
$B_r$. Then, for $r$ small enough,
the intersections $M_i\cap S_r$ are transversal, and moreover, since 
$Q$ is at most countable, we can find $r>0$ as small as we wish with 
 $\gamma_1^r\cap \gamma_2^r=\emptyset $, where 
$\gamma_i^r=M_i\cap S_r$. Moreover, $\gamma_i^r$ will be  transversal to the 
distribution $F_r$ 
of $J$-complex tangent planes to $S_r$, see Lemma A2.1.1 in [I-S]. Let 
$l(\gamma_1^r,\gamma_2^r)$ denote the linking number of the curves 
$\gamma_1^r$ and $\gamma_2^r$. A standard lemma 
from topology (see e.g. [Rf]) says now that 

$$
l(\gamma_1^r,\gamma_2^r)=\sum_{p\in Q\cap B_r}\delta_p.\eqno(4.1)
$$

The positivity of intersections of $J$-complex curves implies that if $M_1$ and 
$M_2$ do intersect at zero, then this linking number is not zero (and is 
positive).

\state Proof of Corollary 1. {\rm All we need to check is that $Y=X\setminus 
D$ does not admit a limiting $J$-complex line. Let $u:\cc \to D$ be a limiting 
$J$-complex line for $Y=X\setminus D$ and let $u(\cc )\subset D_1$ for 
definitivity. $u(\cc )$ cannot be contained 
in  $D_1\setminus Sing(D)$, because this curve is assumed to be hyperbolic. 
Therefore there exists $z_0\in \cc $ such 
that $p:=u(z_0)\in D_1\cap Sing(D)$. Take as $M_1$ the piece of $D_1$ which 
 containing the image of the germ 
of $u$ at $z_0$ and as $M_2$ a germ of any other branch of some $D_i$ passing
through  $p$. 

Chose a sphere $S_r$ as above. Take some $R$ such that $\Delta (R)\supset 
u^{-1}(M_1)$. 
Let $u_n:\Delta (R)\to Y$  be 
$J$-holomorphic mapping approaching uniformly to $u$. Then for $n>>1$ 
the linking number $l(\gamma_n,\gamma_2^r)$ is the same as   
$l(\gamma_1^r,\gamma_2^r)$. Here $\gamma_n$ is an appropriate component of 
$u_n(\Delta (R))\cap S_r$. The positivity of $l(\gamma_1^r,\gamma_2^r)$ implies 
now the positivity of $l(\gamma_n,\gamma_2^r)$ for $n$ large enough. This 
contradicts  the definition of a limiting line, \ie, $u_n(\cc )$ should 
not intersect $D$.

}
\smallskip
\qed

\bigskip\bigskip\bigskip\break

\state Proof of Corollary 2.  {\rm This is rather similar to the previous one.
Suppose we can find a sequence $(\{ D_i^{(k)}\}_{i=1}^5,J_k)$ in $\calm_{\omega
,5l}$, such that $Y_k=(\cc\pp^2\setminus D^{(k)},J_k)$ are  not hyperbolic 
and converge  to $(\{ D_i\}_{i=1}^5,J)$ in $\calm_{\omega
,5l}$, such that $Y=(\cc\pp^2\setminus D,J)$ is hyperbolic. Take 
$J_k$ and $J$-Hermitian metrics $\omega (\cdot ,J_k\cdot )$ and  
$\omega (\cdot ,J\cdot )$ on $\cc\pp^2$. By Corollary 1 there exist
$J_k$-complex line $u_k:\cc \to Y_k$. Gradients of $u_k$ are uniformly 
bounded and therefore some subsequence, stil denoted by $u_k$, converges 
uniformly on compacts in $\cc $. All topologies on the space of holomorphic 
mappings are equivalent (see e.g. Corollary 3.2.2 from [I-S]) and therefore 
$u_k$ converge to some $J$-holomorphic $u:\cc \to \cc\pp^2$ in the $C^1$-sense. 
Since $\Vert du(0)\Vert_h=\lim_{k\to \infty }\Vert du_k(0)\Vert_{h_k}={1\over 2}
$, the map $u$ is a $J$-complex line. $u(\cc )$ cannot be contained 
in $Y$ because of the hyperbolicity of the latter. So $u(\cc )\cap D\not= 
\emptyset $.

Considerations similar to that from the proof of  Corollary 1 show that 
in this case $u_k(\cc )\cap D^{(k)}\not= \emptyset $ for $k>>1$. 
Contradiction.

\smallskip
\qed   
}

\smallskip A similar proof yields one more corollary. Denote by $\calm_{\omega ,
d}$ the Banach manifold of pairs $(J,D)$, where $J$ is an $\omega $-tamed 
a.-c. structure on $\cc\pp^2$ and $D$ is some {\it irreducible } $J$-complex 
curve of degree $d$ in $\cc\pp^2$. Denote by ${\Cal H}_{\omega ,d}$  the subset 
of $(J,D)\in \calm_{\omega ,d}$ such that $(\cc\pp^2\setminus D,J)$ is 
hyperbolically imbedded into $(\cc\pp^2,J)$.

\state Corollary 4.1. {\it The set ${\Cal H}_{\omega ,d}$ is a nonempty open
subset of  $\calm_{\omega ,d}$, provided $d\ge 5$.
}
\smallskip ${\Cal H}_{\omega ,d}$ is not empty by the result of 
Zaidenberg [Za-2], which asserts that there exists an open set in the manifold of 
algebraic curves of degree $d\ge 5$ in $\cc\pp^2$ such that 
$\cc\pp^2\setminus D$ is hyperbolically imbedded into $\cc\pp^2$ for every 
curve $D$ from this set .

\spaceskip=4pt plus3.5pt minus 1.5pt
\font\csc=cmcsc10

\newdimen\length
\newdimen\lleftskip
\lleftskip=3.8\parindent
\length=\hsize \advance\length-\lleftskip
\def\entry#1#2#3#4\par{\parshape=2  0pt  \hsize%
\lleftskip \length%
\noindent\hbox to \lleftskip%
{\bf[#1]\hfill}
{\csc{#2}} 
{\sl{#3}} 
#4 \medskip
}
\ifx \twelvebf\undefined \font\twelvebf=cmbx12\fi

\bigskip\bigskip
\centerline{\bf R E F E R E N C E S}

\nobreak
\bigskip


\entry{Bt}{T.J. Barth.:}{The Kobayashi distance induces the standard topology.}
Proc. Amer. Math. Soc. {\bf 35}, 439-440 (1972).

\entry{Br}{R. Brody.:}{Compact manifolds and hyperbolicity.} Trans. Amer. Math. 
Soc. {\bf 235} 213-219 (1978).

\entry{De}{R. Debalme.:}{Kobayashi hyperbolicity of almost complex manifolds.} 
Preprint of the university of Lille, IRMA {\bf 50} (1999), math.CV/9805130.

\entry{Gr}{M. Green.:}{The hyperbolicity of the complement of 2n+1 hyperplanes in 
general position in $\Bbb P^n$, and related results.} Proc. Amer. Math. Soc.  
{\bf66} 109-113 (1977).

\entry{G}{M. Gromov.:}{Pseudo holomorphic curves in symplectic manifolds.}
Invent. math. {\bf82} 307-347 (1985).

\entry{I-S}{S. Ivashkovich, V. Shevchishin.:}{Complex Curves in Almost-Complex 
Manifolds and Meromorphic Hulls.} Lecture Notes in Schriftenreihe des 
Graduiertenkollegs, Geometrie und Math. Physik, Ruhr-Uni-Bochum, Heft 36, 
1-186 (1999), see also math.CV/9912046.

\bigskip\bigskip\bigskip\break

\entry{Ki}{P. Kiernan.:}{Hyperbolically Imbedded Spaces and the Big Picard 
Theorem.} Math. Ann. {\bf204} 203-209 (1973).

\entry{Ko-1}{S. Kobayashi.:}{Invariant distances on complex manifolds and 
holomorphic mappings.} J. Math. Soc. Japan {\bf19} 460-480 (1967). 

\entry{Ko-2}{S. Kobayashi.:}{Hyperbolic Complex  Spaces.}
Springer (1998).

\entry{Kr}{B. Kruglikov.:}{On the Kobayashi-Royden pseudonorm for almost complex 
manifolds.} math.dg-ga/9708012. 

\entry{KrOv}{B. Kruglikov and M. Overholt.:}{The Kobayashi pseudodistance on 
almost complex manifolds.} Preprint of the university of Tromsoe {\bf19} (1997)
math.dg-ga/9703005. 

\entry{La}{S. Lang.:}{Introduction to complex hyperbolic spaces.}Springer-Verlag 
New York (1987).

\entry{McD}{D. McDuff.:}{Singularities and positivity of intersections of 
J-holomorphic curves.} Holomorphic curves in symplectic geometry, eds. M. Audin 
and J.Lafontaine, Birkhauser, Basel 191-215 (1994). 

\entry{M-W}{M. Micallef , B. White }{The structure of branch points in 
minimal surfaces and in pseudoholomorphic curves.} Ann. Math. {\bf139} 35-85 
(1994).

\entry{Rf}{Rolfsen, D.:}{Knots and links.}Publish or perish, {\bf N7}, (1976).

\entry{Ro}{H.L. Royden.:}{Remarks on the Kobayashi metric.} Springer LNM  
{\bf 185} 125-137 (1970).

\entry{Sk}{J.-C. Sikorav.:}{Some properties of holomorphic curves in almost 
complex manifolds.} Holomorphic curves in symplectic geometry, eds. M. Audin 
and J.Lafontaine, Birkhauser, Basel 165-189 (1994). 

\entry{Za-1}{M. Zaidenberg.:}{Picard's Theorem and Hyperbolicity.} Siberian 
Math. J., {\bf24} 858-867 (1983).

\entry{Za-2}{M. Zaidenberg.:}{Stability of hyperbolic embeddedness and 
construction of examples.} Math. USSR Sbornik, {\bf63}, N 2, 351-361 (1989).

\end